\documentclass[letter,10pt,leqno, english]{amsart}
        \title[The evenness conjecture in equivariant unitary bordism]{The evenness conjecture in equivariant unitary bordism}
      \author{Bernardo Uribe}
      \address{Departamento de Matem\'aticas y Estad\'istica\\
       Universidad del Norte\\ 
       Km. 5 via Puerto Colombia, Barranquilla, Colombia}
      \email{bjongbloed@uninorte.edu.co}
      \urladdr{https://sites.google.com/site/bernardouribejongbloed/}
\thanks{The author acknowledges and thanks the financial support 
provided by the Max Planck Institute for Mathematics, by COLCIENCIAS through the Fondo Nacional de Financiamiento 
para la Ciencia, la Tecnolog\'ia y la Innovaci\'on and by the 
Alexander Von Humboldt Foundation.}

\usepackage{pdfsync}
\usepackage{hyperref}
\usepackage{calc}
\usepackage{enumerate,amssymb,amscd}
\usepackage{color}
\usepackage[arrow,curve,matrix,tips,2cell]{xy}
  \SelectTips{eu}{10} \UseTips
  \UseAllTwocells
\usepackage{tikz}
\usepackage{pdfcolmk}

\DeclareMathAlphabet{\matheurm}{U}{eur}{m}{n}

\DeclareMathOperator{\Ind}{Ind}

\DeclareMathOperator{\map}{map}

\DeclareMathOperator{\Hom}{\textup{Hom}}
\DeclareMathOperator{\Irr}{\textup{Irr}}

  \newcommand{\IC}{\mathbb{C}}

  \newcommand{\IQ}{\mathbb{Q}}
  \newcommand{\IR}{\mathbb{R}}
  \newcommand{\IS}{\mathbb{S}}

  \newcommand{\IV}{\mathbb{V}}

  \newcommand{\IZ}{\mathbb{Z}}

  \newcommand{\calf}{\mathcal{F}}

  \newcommand{\calp}{\mathcal{P}}

\makeatletter
\newcommand{\colim@}[2]{%
  \vtop{\m@th\ialign{##\cr
    \hfil$#1\operator@font colim$\hfil\cr
    \noalign{\nointerlineskip\kern1.5\ex@}#2\cr
    \noalign{\nointerlineskip\kern-\ex@}\cr}}%
}
\newcommand{\colim}{%
  \mathop{\mathpalette\colim@{\rightarrowfill@\textstyle}}\nmlimits@
}
\makeatother

\newcounter{commentcounter}

\theoremstyle{plain}
\newtheorem{theorem}{Theorem}[section]

\newtheorem*{theorem*}{Theorem}
\newtheorem*{mtheorem*}{Main Theorem}

\theoremstyle{definition}
\newtheorem{definition}[theorem]{Definition}

\theoremstyle{remark}
\newtheorem*{summary*}{Summary}

\makeatletter\let\c@equation=\c@theorem\makeatother

\newcommand{\version}[1] 
{\begin{center} Last edited on #1\\
    Last compiled on \today\\
file name: \jobname
  \end{center}
}

\begin{document}

\begin{abstract}  
The evenness conjecture for the equivariant unitary bordism groups states that these bordism
groups are free modules over the unitary bordism ring on even-dimensional generators.
In this paper we review the cases in which the conjecture is known to hold and we highlight
the properties that permit one to prove the conjecture in these cases.   
 \end{abstract}

\maketitle

\section*{Introduction}

 The $G$-equivariant unitary bordism groups for $G$ a compact Lie group are the bordism groups of
$G$-equivariant tangentially stable almost complex manifolds, also known as $G$-equivariant unitary manifolds. These are closed $G$-manifolds $M$ for which
a stable tangent bundle 
$TM \oplus \underline{\IR}^k$, where $\underline{\IR}^k$ denotes the trivial bundle $\IR^k \times M$ with trivial $G$-action, can be endowed with the structure of a $G$-equivariant complex bundle. Two tangentially
 stable almost complex $G$-structures are identified if after stabilization with further 
$G$-trivial $\underline{\IC}$ summands the structures become $G$-homotopic through complex $G$-structures.
Being unitary is inherited by the fixed points sets. Whenever $H$ is a closed subgroup of $G$ the fixed points $M^H$  are also tangentially stable almost complex, and moreover
a $N_{H}$-tubular neighbourhood around $M^H$ in $M$ possesses a complex $N_{H}$ structure \cite[\S XVIII, Proposition 3.2]{May-book}. 

 For a cofibration of $G$-spaces $Y \to X$, the geometric $G$-equivariant unitary bordism groups $\Omega^G_n(X,Y)$ are the
$G$-bordism classes of $G$-equivariant $n$-dimensional manifolds with map   $(M^n, \partial M^n) \to (X,Y)$.

The $G$-equivariant unitary bordism groups of a point $\Omega_*^G$ become a ring under the cartesian
product of manifolds with the diagonal $G$-action, and therefore a module over the unitary bordism ring $\Omega_*$
where we consider a unitary manifold as a trivial $G$-manifold. 
Milnor and Novikov \cite{Milnor, Novikov}, making use of the Adams spectral sequence, showed that the unitary bordism ring 
is a polynomial ring $\Omega_*= \IZ[x_{2i} \colon i \geq 1]$ with one generator in each even degree. In this work we will be interested
in the $\Omega_*$-module structure of the equivariant unitary bordism groups $\Omega_*^G$.

Explicit calculations carried out by Landweber \cite{Landweber-cyclic} in the cyclic case and by Stong \cite{Stong-complex} in the abelian $p$-group case permitted them to conclude that in these cases the equivariant unitary bordism group $\Omega_*^G$ 
is a free $\Omega_*$-module on even-dimensional generators.
Ossa in \cite{Ossa} generalized this result to any finite abelian group and L\"{o}ffler in \cite{Loffler} and Comeza\~na in \cite[\S XXVIII, Theorem 5.1]{May-book} showed that this also holds whenever $G$ is a compact
abelian Lie group. Explicit calculations done for the dihedral groups $D_{2p}$ with $p$ prime by
 \'Angel, G\'omez and the author \cite{AngelGomezUribe}, for groups of order $pq$ with $p$ and $q$ different primes  by Lazarov \cite{Lazarov} and for groups on which all its Sylow subgroups are cyclic by Rowlett \cite{Rowlett-metacyclic} show that for these groups this phenomenon also occurs.
We believe that this property should hold in the $G$-equivariant unitary bordism groups for any compact Lie group $G$, in the same way that
the coefficients for $G$-equivariant K-theory are trivial in odd degrees and a free module over the integers on even degrees. 

The theme of this work is the
\begin{center}
 {\bf  Evenness conjecture for the equivariant unitary bordism groups}
\end{center}
which states that the $G$-equivariant unitary bordism group is a free $\Omega_*$-module
on even-dimensional generators whenever $G$ is a compact Lie group. Rowlett explicitly mentions this conjecture
in his work of 1980  \cite{Rowlett-metacyclic}  and
 later Comeza\~na  in his work  of 1996 \cite[\S XXVIII.5]{May-book}.
We also believe that this conjecture holds in general and we do hope that this paper will help to spread it
to the mathematical community for its eventual proof.

In this work we survey the original proofs of the known cases of the evenness conjecture 
for finite groups. We start in section \ref{section bordism} with the definition of the equivariant unitary bordism groups for pairs of families and
the long exact sequence associated to them. In section \ref{section bundles} we review the decomposition of 
equivariant complex vector bundles  restricted to fixed point sets done by \'Angel, G\'omez and the author \cite{AngelGomezUribe} and how this decomposition
allows one to write the equivariant unitary bordism groups of adjacent pair of families as the bordism groups
of equivariant classifying spaces. In section \ref{section abelian} we review the proofs of the evenness conjecture 
done by Landweber \cite{Landweber-cyclic} for cyclic groups, by Stong \cite{Stong-complex} for abelian $p$-groups and by Ossa \cite{Ossa} for general finite abelian groups. In section \ref{section non abelian} we review the proof of the
evenness conjecture for groups for which all its Sylow subgroups are cyclic done by Rowlett \cite{Rowlett-metacyclic} and we finish in section \ref{section conclusion} with some conclusions.

We would like to thank Prof. Peter Landweber for reading this manuscript and for providing valuable comments and suggestions
that improved the presentation of this survey.

\section{Equivariant unitary bordism for families of subgroups} \label{section bordism}

 To study the equivariant bordism groups Conner and Floyd introduced
the study of bordism groups of manifolds with prescribed isotropy groups \cite[\S 5]{ConnerFloyd-Odd}.
A family of subgroups  $\calf$ of $G$ is a set of subgroups of $G$ which is closed under taking subgroups and under conjugation.
The classifying space for the family $E\calf$ is a $G$-space which is terminal in the category of $\calf$-numerable
$G$-spaces \cite[\S 1, Theorem 6.6]{tomDieck-transformation} and characterized by the following properties on fixed point sets:
$E\calf^H\simeq*$ if $H \in \calf$ and $E\calf^H = \emptyset$ if $H \notin \calf$.  This classifying space may
be constructed in such a way that whenever
 $\calf' \subset \calf$, the induced map $E \calf' \to E\calf$ is a
 $G$-cofibration.

The equivariant unitary bordism groups for pairs of families may be defined as follows
$$\Omega_*^G[\calf,\calf'](X,A):= \Omega_*^G(X \times E\calf,X \times E\calf' \cup A \times E\calf'),$$
see \cite[page 310]{tomDieck-Orbit-I}, or alternatively they
may be defined in a geometric way as in \cite[\S 2]{Stong-complex}. 

A $(\calf, \calf')$ {\it free geometric unitary bordism element} of $(X,A)$ is an equivalence class of 4-tuples $(M,M_0,M_1,f)$, where:
\begin{itemize}
\item $M$ is an $n$-dimensional $G$-manifold endowed with tangentially stable almost $G$ structure which is moreover $\calf$-free, i.e. such that all isotropy groups
$G_m=\{ g \in G \ | \ gm=m \}$ for $m \in M$ belong to $\calf$, and such that $f:M \to X$ is $G$-equivariant; and
\item $M_0, M_1$ are compact submanifolds of the boundary of $M$, with $\partial M = M_0 \cup M_1$, $M_0 \cap M_1 = \partial M_0=\partial M_1$
having tangentially stable almost complex structures induced from $M$, both $G$-invariant, such that $f(M_1) \subset A$ and $M_0$ if $\calf'$-free,
i.e. all isotropy groups of $M_0$ belong to $\calf'$.
\end{itemize}

Two four-tuples $(M,M_0,M_1,f)$ and $(M',M'_0,M'_1,f)$ are equivalent if there is a 5-tuple $(V,V^+,V_0,V_1,F)$ where
\begin{itemize}
\item $V$ is a $\calf$-free manifold and $F : V \to X$ is a $G$-equivariant map;
\item The boundary of $V$ is the union of $M$, $M'$ and $V^+$ with $M \cap V^+=\partial M$, $M' \cap V^+= \partial M'$, $M \cap M' = \emptyset$,
 $V^+ \cap(M \cup M') = \partial V^+$, with $V$ inducing the tangentially stable almost complex $G$-structure on $M$ and the opposite one on $M'$;
  $V^+$ is 
$G$-invariant and $F$ restricts to $f$ in $M$ and to $f'$ on $M'$; and
\item $V^+$ is the union of the $G$-invariant submanifolds $V_0$, $V_1$ with intersection a submanifold $V^-$ in their boundaries, such that
$\partial V_i=M_i \cup V^-\cup M_i'$, $M_i \cap V^-= \partial M_i$,  $M'_i \cap V^-= \partial M'_i$ with $V_0$ is $\calf'$-free and $F(V_1) \subset A$.
\end{itemize}

\begin{definition}
The set of equivalence classes of $n$-dimensional $(\calf, \calf')$-free geometric unitary bordism elements of $(X,A)$, consisting of
classes $(M,M_0,M_1,f)$ where the dimension of $M$ is $n$, and under the operation of disjoint union, forms an abelian
group denoted by $$\Omega^G_n \{\calf ,\calf' \}(X,A).$$ Call these groups the geometric $G$-equivariant unitary bordism groups of the pair $(X,A)$
restricted to the pair of families $\calf' \subset \calf$.
\end{definition}
 Note that if $N$ is a tangentially stable almost complex closed manifold, we can
define $N \cdot (M,M_0,M_1,f) = (N \times M,N \times M_0,N \times M_1,f \circ \pi_M)$ thus making  $\Omega^G_* \{\calf ,\calf' \}(X,A)$
a graded module over the unitary bordism ring $\Omega_*$.

The covariant functor $\Omega^G_* \{\calf ,\calf' \}$ defines a $G$-equivariant homology theory
 \cite[Proposition 2.1]{Stong-complex}, the boundary map on $A$
\begin{align*}
\delta : \Omega^G_n \{\calf ,\calf' \}(X,A) & \to \Omega^G_{n-1} \{\calf ,\calf' \}(A,\emptyset) \\
(M,M_0,M_1,f) & \mapsto (M_1,\partial M_1,\emptyset,f|_{M_1})
\end{align*}
defines the long exact sequence in homology for pairs
$$\cdots \Omega^G_n \{\calf ,\calf' \}(X,A) \stackrel{\delta}{\to} \Omega^G_{n-1} \{\calf ,\calf' \}(A,\emptyset) \to \Omega^G_{n-1} \{\calf ,\calf' \}(X,\emptyset) \to \cdots,$$
and for families  $\calf'' \subset \calf' \subset \calf$, choosing the boundary which is $\calf'$-free 
\begin{align*}
\partial : \Omega^G_n \{\calf ,\calf' \}(X,A) & \to \Omega^G_{n-1} \{\calf',\calf'' \}(X,A) \\
(M,M_0,M_1,f) & \mapsto (M_0, \emptyset, \partial M_0,f|_{M_0})
\end{align*}
one obtains by \cite[Proposition 2.2]{Stong-complex} the long exact sequence in homology for families
$$\cdots \Omega^G_n \{\calf ,\calf' \}(X,A) \stackrel{\partial}{\to} \Omega^G_{n-1} \{\calf' ,\calf'' \}(X,A) \to \Omega^G_{n-1} \{\calf ,\calf'' \}(X,A) \to \cdots.$$

The bordism condition restricted to the non-relative case $\Omega^G_* \{\calf ,\calf' \}(X)$ can be read as the set of bordism classes of maps $f:M \to X$ such
that $M$ is $\calf$-free and $\partial M$ is $\calf'$-free with $M$ endowed with a tangentially stable almost complex $G$-structure. Two become  equivalent if there exists a $G$-manifold $F:V \to X$ which is $\calf$-free such that $\partial V=M \cup M' \cup V^+$
and $M \cap V^+=\partial M$, $M' \cap V^+= \partial M'$, $M \cap M' = \emptyset$, $V^+ \cap(M \cup M') = \partial V^+$, with
the property that $F$ restricts to $f$ on $M$ and to $f'$ on $M'$ and with $V^+$ $\calf'$-free.

In \cite[Satz 3]{tomDieck-Orbit-I} it is shown that the canonical map that one can define 
\begin{align} \label{iso geometric families EF}
\mu : \Omega^G_n \{\calf ,\calf' \}(X,A) \to \Omega^G_n [\calf ,\calf' ](X,A)
\end{align}
becomes  a natural isomorphism of  homology theories.

A key fact about the $(\calf, \calf')$-free geometric unitary bordism elements of $X$ is the following result 
proven in \cite[Lemma 5.2]{ConnerFloyd-Odd}. Whenever $(M^n,\partial M^n ,f)$ is a $(\calf, \calf')$-free geometric unitary bordism element of $X$ and $W^n$ a compact
manifold with boundary regularly embedded in the interior of $M^n$ and invariant under the $G$-action, such that $G_m \in \calf'$
for all $m \in M^n \backslash W^n$, then $[M^n,\partial M^n ,f]=[W^n, \partial W_n, f|_{W^n}]$ in $\Omega^G_n \{\calf ,\calf' \}(X)$.

Whenever the pair of families $\calf' \subset \calf$ differ by  a fixed group $A$, i.e.
$\calf \backslash \calf' = (A)$ with $(A)$ the set of subgroups of $G$ conjugate to $A$, then the pair
$(\calf, \calf')$ is called an {\it adjacent pair of families of groups}. In the case that $A$ is normal in $G$
a  $(\calf, \calf')$-free geometric unitary bordism class $[M, \partial M, f]$ of $X$ is equivalent 
to $\sum_{j=1}^l [U_j, \partial U_j, f|_{U_j}]$ where the $U_j$'s are disjoint $G$-equivariant tubular neighborhoods
of the $M^A_j$'s and these sets are the connected components of the $A$-fixed point set $M^A$.
Since the normal bundle of the fixed point set $M^A_j$ may be classified by a map
to an appropriate classifying space, the groups $\Omega^G_* \{\calf, \calf'\}(X)$ become
isomorphic to the direct sum of $G/A$-free equivariant unitary bordism groups of the product of $X^A$ with an appropriate
classifying space (see \cite[Theorem 4.5]{AngelGomezUribe}). To introduce this result
we need to understand how the fixed points of universal equivariant bundles behave. This is the subject of the next section.

\section{Equivariant vector bundles and fixed points} \label{section bundles}
 
\subsection{Complex representations} Let $G$ be a compact Lie group and $A$ a closed and normal subgroup of $G$ fitting into the exact sequence
$$1 \to A \to G \to Q \to 1.$$
Let $\rho: A \to U(V_\rho)$ be an irreducible unitary representation of $A$, denote by $\Irr(A)$ the set of isomorphism classes of irreducible representations of $A$
 and let $W$ be a finite dimensional complex $G$-representation.
 Then we have an isomorphism of $A$-representations
$$\bigoplus_{\rho \in \Irr(A)} V_\rho \otimes \Hom_A(V_\rho,W) \stackrel{\cong}{\to} W.$$

The group $G$ acts on the set of $A$-representations
$$(g\cdot \rho) (a) := \rho(g^{-1}ag)$$
 and therefore it acts on $\Irr(A)$. Denote by 
 $G_\rho := \{g \in G \ | \ g \cdot \rho \cong \rho \}$ the isotropy group of the isomorphism class of $\rho$
 and denote 
  $Q_\rho:= G_\rho/A$.
 
 If $g \cdot \rho \cong \rho$ then there exists $ M\in U(V_\rho)$ such that $g \cdot \rho (a) = M^{-1}\rho(a) M$.
 Since this matrix $M$ is unique up to a central element, we obtain a homomorphism $f: G_\rho \to PU(V_{\rho})$
 which fits into the following diagram 
\[
\xymatrixrowsep{0.5cm}
\xymatrix{ A \ar[d]_\rho \ar@{^{(}->}[r]^{\iota} & G_\rho \ar[d]^f \\
U(V_\rho) \ar[r]^p & PU(V_{\rho}),
 }
\] 
thus making $V_\rho$ into a projective $G_\rho$-representation.

Define the $\IS^1$-central extension $ \widetilde{G}_{\rho} : = f^*U(V_\rho)$ of $G_\rho$ which fits into the
following diagram
$$
\xymatrixrowsep{0.4cm}
\xymatrix{
& \IS^1  \ar[d]  &  \IS^1 \ar[d] \\ 
A \ar[r]^{\widetilde{\iota}}  \ar[d]^= & \widetilde{G}_{\rho}
\ar[r]^{\widetilde{f}} \ar[d] & U(V_\rho) \ar[d] \\
A \ar[r]^\iota  & G_\rho \ar[r]^f & PU(V_\rho),
}
$$
endowing $V_\rho$ with the structure of a $\widetilde{G}_{\rho}$-representation where $\IS^1$ acts by multiplication with scalars.

The vector space  $\Hom_A(V_\rho, W)$
is also a $\widetilde{G}_{\rho}$-representation where
for $\phi \in  \Hom_A(V_\rho, W)$ and $\widetilde{g} \in \widetilde{G}_{\rho}$ we set
$$( \widetilde{g} \bullet \phi)(v) := g \phi(\widetilde{f}(\widetilde{g})^{-1} v).$$
It follows that $A$ acts trivially on $\Hom_A(V_\rho, W)$ and moreover the elements of $\IS^1$
 act by multiplication with their inverse.
 
 Hence $V_\rho \otimes \Hom_A(V_\rho, W)$
is a $G_\rho$ representation, where  $\Hom_A(V_\rho, W)$ is a $\widetilde{Q}_{\rho}:= \widetilde{G}_{\rho}/A$ representation where $\IS^1$
acts by multiplication of the inverse. Here $\widetilde{Q}_{\rho}$
 is an $\IS^1$-central extension of $Q_\rho$.

Since the isotropy group $G_\rho$ contains the connected component of the identity in $G$,
the index $[G \colon G_\rho]$ is finite and we may induce the $G_\rho$-representation
$V_\rho \otimes \Hom_A(V_\rho, W)$ to $G$ thus obtaining the following result.
\begin{theorem} There is 
a canonical
isomorphism of $G$-representations
$$ \bigoplus_{\rho \in G \backslash \Irr(A)} \Ind_{G_\rho}^G \left( V_\rho \otimes \Hom_A(V_\rho, W) \right) \cong W$$
where $\rho$ runs over  representatives of the orbits of the action of $G$ on $\Irr(A)$. 
\end{theorem}

\subsection{Equivariant complex bundles} The previous result generalizes to equivariant complex vector bundles, but prior to showing this generalization we need
to recall the multiplicative induction map introduced in \cite[\S 4]{Bix-tomDieck}. Let $H$ be a closed
subgroup of the compact Lie group $G$. The right adjoint to the restriction functor $r_H^G$ from
$G$-spaces to $H$-spaces is called the multiplicative induction functor and takes an $H$-space
$Y$ and returns the $G$-space
$$m_H^G(Y):=\map(G,Y)^H$$
of $H$-equivariant maps from $G$ to $Y$, with $G$ considered as an $H$-space by left multiplication. The $G$-action
on $m_H^G(Y)$ is given by $(g \cdot f)(k) := f(kg)$, $m_H^G(Y)$ is homeomorphic to the space
of sections of the projection map $G \times_H Y \to G/H$ and, in the case that $G/H$ is finite,
it is homeomorphic to $[G:H]$ copies of $Y$.

  There is a homeomorphism
$$\map(X,m_H^G(Y))^G \stackrel{\cong}{\to} \map(r_H^G(X), Y)^H, \ \ F \mapsto (x \mapsto F(x)(1_G))$$
whose inverse maps $f$ to $m_H^G(f) \circ p_H^G$ where $p_H^G:X \to m_H^G(r_H^G(X))$, $p_H^G(x)(g)=gx$,
is the unit of the adjunction.

Now consider a $G$-space $X$ on which the closed and normal subgroup $A$ acts trivially. Take a $G$-equivariant
complex vector bundle $p: E \to X$ and assume that $E$ has an hermitian metric in such a way that
$G$ acts through unitary matrices on the complex fibers. For a complex $A$-representation
$\rho : A \to U(V_\rho)$ denote by $\IV_{\tau}$ the trivial $A$-vector bundle $\pi_2 : V_\rho \times X \to X$.

The complex vector bundle $\Hom_{A}(\IV_{\rho},E)$ is a  $\widetilde{Q}_{\rho}$-equivariant complex vector bundle where $\IS^1$ acts on the fibers by multiplication of the inverse,
 $\IV_{\rho} \otimes\Hom_{A}(\IV_{\rho},E)$ is a $G_\rho$-equivariant complex vector bundle and
$$(p_{G_\rho}^G)^* \left(m_{G_\rho}^G(\IV_{\rho} \otimes\Hom_{A}(\IV_{\rho},E)) \right) \to X$$
is a $G$-equivariant complex vector bundle over $X$. 

\begin{theorem}\cite[Theorem 2.7]{AngelGomezUribe}
\label{thm decomposition} Let $G$ be a compact Lie group, $A$ a closed and
normal subgroup, $X$ a $G$-space on which $A$ acts trivially and $E \to X$ a $G$-equivariant complex
vector bundle. Then there is an isomorphism of $G$-equivariant complex vector bundles
$$\bigoplus_{\rho \in G \backslash \Irr(A)} (p_{G_\rho}^G)^* \left(m_{G_\rho}^G(\IV_{\rho} 
\otimes\Hom_{A}(\IV_{\rho},E)) \right)\stackrel{\cong}{\to} E$$
where $\rho$ runs over representatives of the orbits of the $G$-action on the 
set of isomorphism classes of $A$-irreducible representations.
\end{theorem}

With the same hypothesis as in the previous theorem, there is an induced decomposition in equivariant K-theory
$$ K^*_G(X) \cong \bigoplus_{\rho \in G \backslash \Irr(A)}
{}^{\widetilde{Q}_{\rho}} K^*_{Q_\rho}(X), \ \ E \mapsto \bigoplus_{\rho \in G \backslash \Irr(A)}
\Hom_{A}(\IV_{\rho},E)$$
where ${}^{\widetilde{Q}_{\rho}} K^*_{Q_\rho}(X)$ is the $\widetilde{Q}_{\rho}$-twisted $Q_\rho$-
equivariant K-theory of $X$ which is built out of the Grothendieck group of $\widetilde{Q}_{\rho}$-equivariant
complex vector bundles over $X$ on which the central $\IS^1$ acts on the fibers by multiplication.

\subsection{Classifying spaces} The decomposition described above can also be written at the level of classifying spaces; let us
set up the notation first. Let $G$ be a compact Lie group and $\widetilde{G}$ a $\IS^1$-central group extension
of $G$. Let $\widetilde{\bf{C}}^\infty$ be a countable direct sum of all complex irreducible $\widetilde{G}$ representations
on which elements of $\IS^1$ act by multiplication with their inverse. Denote by  ${}^{\widetilde{G}}B_GU(n)$ the Grassmannian
of $n$-planes of  $\widetilde{\bf{C}}^\infty$ and denote by ${}^{\widetilde{G}}\gamma_GU(n)$ the canonical
$n$-plane bundle over ${}^{\widetilde{G}}B_GU(n)$.
The complex vector bundle
$$ \IC^{n} \to
{}^{\widetilde{G}}\gamma_GU(n) \to 
{}^{\widetilde{G}}B_GU(n)$$
is a universal $\widetilde{G}$-twisted $G$-equivariant complex vector bundle of rank $n$.
Denote by $\gamma_GU(n) \to B_GU(n)$ the universal $G$-equivariant complex vector bundle of rank $n$.

For a closed subgroup $A$ of $G$, let $N_A$ denote the normalizer of $A$ in $G$ and $W_A:= N_A/A$.
Consider the fixed point set $B_GU(n)^A$ and the restriction $\gamma_GU(n)|_{B_GU(n)^A}$
of the universal bundle to this fixed point set. For $\rho \in \Irr(A)$, by the arguments above
we have that
$$\Hom_{A}(\IV_{\rho},\gamma_GU(n)|_{B_GU(n)^A})$$
 is a $(\widetilde{W_A})_\rho$-twisted
$(W_A)_\rho$-equivariant complex bundle, but since the space $B_GU(n)^A$ is not necessarily connected,
it may not have constant rank. Therefore Theorem \ref{thm decomposition} implies the following
equivariant homotopy equivalence.

\begin{theorem} \cite[Theorems 3.3 \& 3.5]{AngelGomezUribe}
\label{thm homotopy decomposition}
There are $W_A$-equivariant homotopy equivalences
$$\bigsqcup_{n=0}^{\infty}\gamma_GU(n)^A \simeq \left(
\bigsqcup_{n=0}^{\infty}\gamma_{W_A}U(n_1) \right)\times
\prod_{\rho \in W_A \backslash \Irr(A) \atop \rho \ne1} m^{W_A}_{(W_A)_\rho} 
\left( \bigsqcup_{n_\rho=0}^{\infty}{}^{(\widetilde{W_A})_\rho}B_{(W_A)_\rho}U(n_\rho) \right),$$
$$\bigsqcup_{n=0}^{\infty}B_GU(n)^A \simeq 
\prod_{\rho \in W_A \backslash \Irr(A)} m^{W_A}_{(W_A)_\rho} 
\left( \bigsqcup_{n_\rho=0}^{\infty}{}^{(\widetilde{W_A})_\rho}B_{(W_A)_\rho}U(n_\rho) \right).$$
\end{theorem}

If $G$ is abelian then $A$ is normal, $G$ acts trivially on $\Irr(A)$ and all the irreducible
representations are 1-dimensional. Therefore we get a $G/A$-homotopy equivalence
$$\gamma_GU(n)^A \simeq \bigsqcup_{(n_\rho)_{\rho \in \Irr(A)} \atop \sum_\rho n_\rho = n} \left(
\gamma_{G/A}U(n_1) \times \prod_{\rho \in \Irr(A) \atop \rho \neq 1}
B_{G/A}U(n_\rho)\right).$$
In order to get a similar formula for the case in which $G$ is not abelian we need to introduce further notation
and make some choices. Let 
$\calp(n,A)$ be the set of arrangements of non-negative integers $(n_\rho)_{\rho \in  \Irr(A)}$ such that 
$$\sum_{\rho  \in \Irr(A)} n_\rho |\rho|=n,$$ 
then non-equivariantly there is a homotopy equivalence
   $$B_GU(n)^A \simeq \bigsqcup_{{(n_\rho) \in \calp(n,A)   }}  \prod_{\rho \in  \Irr(A)} 
\left({}^{(\widetilde{W_A})_\rho}B_{(W_A)_\rho}U(n_\rho) \right).$$
The group $W_A$ acts on $\calp(n,A)$ on the right by permuting the arrangements, i.e. the action of $g \in W_A$
on the arrangement $(n_\rho)$ is the arrangement
  $(n_\rho) \cdot g:= (n_{g\cdot \rho})$ meaning that it has the number $n_{g \cdot \rho}$ in the coordinate  $\rho$.
 Denote by $(W_A)_{(n_\rho)}$ the isotropy group of the arrangement $(n_\rho)$. Rearranging the terms
 we obtain the following $W_A$-equivariant homotopy equivalence
 \begin{align} 
 &B_GU(n)^A \simeq  \label{homotopy type BGU(n)A} \\  \nonumber
  \bigsqcup_{{(n_\rho) \in \calp(n,A) / W_A   }} &
  W_A \underset{(W_A)_{(n_\rho)}}{\times} \left(
   \prod_{\rho \in (W_A)_{(n_\rho)} \backslash \Irr(A)} m_{(W_A)_\rho \cap (W_A)_{(n_\rho)}}^{(W_A)_{(n_\rho)}}
\left({}^{(\widetilde{W_A})_\rho}B_{(W_A)_\rho}U(n_\rho) \right) \right)
 \end{align}
where $(n_\rho)$ runs over representatives of the orbits of the action of $W_A$  on $\calp(n,A)$,
and $\rho$ runs over representatives of the orbits of the action of $(W_A)_{(n_\rho)}$ on $\Irr(A)$.
  
For the calculation of the equivariant unitary bordism of adjacent families of groups we need
to consider only the arrangements of non-negative integers $(n_\rho)$ such that the number associated
to the trivial representation is zero, i.e. $n_1=0$. Denote by $\overline{\calp}(n,A)$ the set of
arrangements $(n_\rho)$ such that $n_1=0$ and define the $W_A$-space:
\begin{align}
 &C_{N_A,A}(k) :=  \label{space CGA(k)} \\ \nonumber
  \bigsqcup_{{(n_\rho) \in \overline{\calp}(k,A) / W_A   }}&
  W_A \underset{(W_A)_{(n_\rho)}}{\times} \left(
   \prod_{\rho \in (W_A)_{(n_\rho)} \backslash \Irr(A)  \atop \rho \neq 1} m_{(W_A)_\rho \cap (W_A)_{(n_\rho)}}^{(W_A)_{(n_\rho)}}
\left({}^{(\widetilde{W_A})_\rho}B_{(W_A)_\rho}U(n_\rho) \right) \right)
 \end{align}

Therefore we have the following $W_A$-homotopy equivalence
\begin{align} \label{formula with C}
\gamma_GU(n)^A \simeq
\bigsqcup_{k=0}^n \gamma_{W_A}U(n-k) \times C_{N_A,A}(k)
\end{align}
such that in the case that $G$ is abelian we have the simple formula
\begin{align}
C_{G,A}(k) = \bigsqcup_{{(n_\rho) \in \overline{\calp}(k,A)  }} 
 \prod_{\rho \in \Irr(A) \atop \rho \neq 1}
 B_{G/A}U(n_\rho).
 \end{align}

 Now we are ready to state the relation between the $G$-equivariant unitary bordism groups
of adjacent pair of families of groups and the classifying spaces defined above.

\begin{theorem} \cite[Corollary 4.6]{AngelGomezUribe} \label{theorem adjacent}
Let $G$ be a finite group, $X$ a $G$-space and $(\calf, \calf')$ an adjacent pair of families differing
by the conjugacy class of the subgroup $A$, then there is an isomorphism 
$$\Omega^G_n\{\calf, \calf'\}(X) \cong 
\bigoplus_{0 \leq 2k \leq n} \Omega^{W_A}_{n-2k}\{\{1\}\}(X^A \times  C_{N_A,A}(k) )$$
where $\{1\}$ is the family of subgroups of $W_A$ which only contains the trivial group.
\end{theorem}
Take a bordism class $[M, \partial M , f:M \to X]$ in $\Omega^G_n\{\calf, \calf'\}(X)$ and note
that $M^A \cap M^{gAg^{-1}} = \emptyset$ whenever $g$ does not belong to $N_A$. Then choose
an $N_A$-equivariant tubular neighbourhood $U$ of $M^A$ such that its $G$-orbit $G \cdot U$ is
a $G$-equivariant tubular neighbourhood of $G \cdot M^A$ and such that 
$$G \underset{N_A}{\times} U \stackrel{\cong}{\to} G \cdot U, \ \ [(g,u)] \to gu$$
is a $G$-equivariant diffeomorphism. The assignment $[M, \partial M , f:M \to X] \mapsto
[U, \partial U , f|_U:U \to X]$ induces an isomorphism
$$\Omega^G_n\{\calf, \calf'\}(X) \stackrel{\cong}{\to} \Omega^{N_A}_n\{\calf|_{N_A}, \calf'|_{N_A}\}(X).$$

Let $M^A_{n-2k}$ denote the component of $M^A$ which is a $(n-2k)$-dimensional $W_A$-free manifold and
such that $M^A= \bigcup_{0 \leq 2k \leq n} M^A_{n-2k}$. The tubular neighbourhood $U$
is $N_A$-equivariantly diffeomorphic to $ \bigcup_{0 \leq 2k \leq n} D(\nu_{n-2k})$ where
$\nu_{n-2k} \to M^A_{n-2k}$ is the normal bundle of the inclusion $M^A_{n-2k} \to M$. Since
the trivial $A$-representation does not appear on the fibers of the normal bundles, by Theorem
\ref{thm homotopy decomposition} and formula \eqref{formula with C} we know that
the bundle $\nu_{n-2k}$ is classified by a $W_A$-equivariant map $h_{n-2k}:M^A_{n-2k} \to C_{N_A,A}(k)$.
The bordism class $[ M^A_{n-2k}, f|_{M^A_{n-2k}} \times h_{n-2k} : M^A_{n-2k} \to X^A \times  C_{N_A,A}(k)]$
belongs to $\Omega^{W_A}_{n-2k}\{\{1\}\}(X^A \times  C_{N_A,A}(k) )$ and the assignment
$$[M, \partial M , f:M \to X] \mapsto \bigoplus_{0 \leq 2k \leq n} [ M^A_{n-2k}, f|_{M^A_{n-2k}}\times h_{n-2k} : M^A_{n-2k} \to X^A \times  C_{N_A,A}(k)]$$
induces the desired isomorphism.

\section{The evenness conjecture for finite abelian groups} \label{section abelian}

In this section we will outline the main ingredients used by Landweber \cite{Landweber-cyclic} in the cyclic group case,
Stong \cite{Stong-complex} in the $p$-group case and Ossa \cite{Ossa} in the general case to show that the evenness conjecture
holds for finite abelian groups. The conjecture also holds for compact abelian groups, L\"{o}ffler \cite{Loffler} showed it for
the homotopic $G$-equivariant unitary bordism groups in the case that $G$ is a unitary torus, and Comeza\~na \cite[\S XXVIII]{May-book} generalized it to any compact abelian group. Comeza\~na furthermore showed that the map
from the $G$-equivariant unitary bordism groups to the homotopic ones is injective whenever $G$ is compact abelian thus proving
the evenness conjecture for any compact abelian group.
In this work we will address the finite group case. 

Prior to addressing the study of the $G$-equivariant unitary bordism groups for finite abelian groups we need to 
recall some results on the unitary bordism groups. 

Thom's remarkable Theorem \cite{Thom} shows that the unitary bordism groups $\Omega_*$ can
be calculated as the stable homotopy groups $\lim_k\pi_{n+k}(MU(k))$ of the Thom spaces $MU(k)$ of the canonical
complex vector bundles over $BU(k)$. Milnor in \cite[Theorem 3]{Milnor} 
showed that these stable homotopy groups are zero if $n$ is odd and free abelian if $n$ is even with
a number of generators equal to the number of partitions of $n/2$. Independently Novikov \cite[Theorem 1]{Novikov1} showed that as a ring the unitary bordism groups are isomorphic to the ring of polynomials over the integers with generators
$x_{2i}$ of degree $2i$ for $i\geq 1$. The spectrum $MU$ that the Thom spaces $MU(k)$ define permitted
Atiyah \cite{Atiyah-bordism} to define the homotopy unitary bordism groups $MU_*(X)$ and the homotopy unitary cobordism 
groups $MU^*(X)$ of a space $X$ as a generalized homology and cohomology theory respectively. Thom's theorem implies that for $X$ a CW-complex the
unitary bordism groups over $X$ are equivalent to the homotopic ones: $\Omega_*(X) \cong MU_*(X)$
via the Thom-Pontrjagin map.

The Atiyah-Hirzebruch spectral sequence \cite{AtiyahHirzebruch} (cf. \cite[\S 4.2]{Kochman}) applied to the unitary bordism groups of a CW-complex $X$ produces a spectral sequence which converges to $\Omega_*(X)$ and whose second page is $E^2_{p,q}\cong H_p(X; \Omega_q)$; let us call this spectral sequence
the {\it bordism spectral sequence}. The Thom homomorphism
$$\mu: \Omega_*(X) \to H_*(X; \IZ), \ \ [M, f:M \to X] \mapsto f_*[M],$$
which takes a unitary bordism element in $X$ and maps it to the image under $f$ of the fundamental
class $[M] \in H_*(M; \IZ)$, is a natural transformation of homology theories and is also
the edge homomorphism $\Omega_*(X) \to E^2_{*,0} \cong H_*(X; \IZ)$ of the spectral sequence. 

Whenever $X$ is a CW-complex whose homology $H_*(X; \IZ)$  is free abelian then by \cite[Lemma 3.1]{ConnerSmith}
the bordism spectral sequence collapses, the unitary bordism group $\Omega_*(X)$ is a free
$\Omega_*$-module and the homomorphism induced by the Thom map
$$\widetilde{\mu}: \IZ {\otimes}_{\, \Omega_*} \Omega_*(X) \to H_*(X; \IZ)$$
is an isomorphism.

Applying the bordism spectral sequence to the unitary bordism groups of $BU(n)$ it is shown in
\cite[Proposition 4.3.3]{Kochman} that $\Omega_*(BU(n))$ is a free $\Omega_*$-module
with basis
$$\Omega_*(BU(n)) \cong \Omega_* \{ \alpha_{k_1}\alpha_{k_2} \dots \alpha_{k_n} \colon k_1 \leq \cdots \leq k_n \}$$
where $\alpha_{k_1}\alpha_{k_2} \dots \alpha_{k_n}$ is the unitary bordism class of the bordism element
$$(\IC P^{k_1} \times \cdots \times \IC P^{k_n} ,F: \IC P^{k_1} \times \cdots \times \IC P^{k_n} \to BU(n))$$
 where the map $F$ classifies the canonical rank $n$ complex vector bundle over the product of projective spaces.
 
 In \cite[Proposition 3.6]{ConnerSmith} it is shown that if $X$ is a finite CW-complex such that the Thom
 homomorphism is surjective then the bordism spectral sequence collapses. Whenever $BG$ is the classifying
 space of a finite group $G$ Landweber showed in \cite[Theorem 3]{Landweber-complex} that the following conditions
 are equivalent:
 \begin{itemize}
 \item The Thom homomorphism $\mu : \Omega_*(BG) \to H_*(BG; \IZ)$ is surjective.
 \item The bordism spectral sequence collapses.
 \item $G$ has periodic cohomology, i.e. every abelian subgroup of $G$ is cyclic.
 \item $H^n(BG; \IZ)=0$ for all odd $n$.
 \item The projective dimension of $\Omega_*(BG)$ as a $\Omega_*$-module is  $1$ or $0$.
  \end{itemize} 

 The previous result implies that whenever we consider the cyclic group  $G= \IZ/k$ of order $k$,
 the bordism classes $[L^{2n+1}(k), \iota: L^{2n+1}(k) \to B \IZ/k]$ of the lens spaces $L^{2n+1}(k):= S^{2n+1}_k/ (\IZ/k)$ ,
where $S^{2n+1}_k$ denotes the sphere of unit vectors in $\IC^{n+1}$ with the $\IZ/k$-action
given by multiplication of the root of unity $e^{\frac{2 \pi i}{k}}$, generate $\Omega_*(B \IZ/k)$ as a $\Omega_*$-module. 
 
 One property of finite abelian groups that will be used is the following. If $A$ is a subgroup
 of a finite abelian group $G$ and $\Gamma$ is a product of classifying spaces of the form $B_GU(k)$, then by Theorem \ref{thm homotopy decomposition} the fixed point set $\Gamma^A$ is a product
 of classifying spaces of the form $B_{G/A}U(l)$. This fact allows one  to use 
 an induction hypothesis when calculating the equivariant unitary bordism groups
 of products of spaces of the form $B_GU(k)$.
 
 Now we can start the proof of the evenness conjecture for finite abelian groups. First
 we will handle the case of cyclic $p$-groups following \cite{Landweber-cyclic}, then we will
 review the case of general abelian $p$-groups following \cite{Stong-complex} and we will
 prove the general case using a simple argument on localization shown in \cite{Ossa}.

 \subsection{Cyclic $p$-groups}\label{subsection cyclic} Let $G$ be a cyclic group of order $p^s$ a power of the prime $p$. Let $\Gamma:= \prod_{i=1}^lB_GU(k_i)$ be a product of spaces of the form
 $B_GU(k)$ and $\calf_t= \{ H \subset G \colon |H| \leq p^t\}$ the family of of subgroups
 or order bounded by $p^t$; the family $\calf_s$ is the family of all subgroups of $G$ and therefore
 $\Omega_*^G( \ )  =  \Omega_*^G\{ \calf_s \}( \ )$. Let us split $\Omega_*^G( \ ) = \Omega_+^G( \ ) \oplus \Omega_-^G( \ )$
 where $ \Omega_+^G( \ )$ denotes the even degree bordism groups
 and $\Omega_-^G( \ )$ the odd degree ones. We will prove by induction
 on the size of the group that   for any $0 \leq t < s$ the following properties hold:
\begin{itemize}
\item 
 $\Omega^G_*\{\calf_s,\calf_t\}(\Gamma)$ is a free $\Omega_*$-module on even-dimensional generators.
 \item $\Omega^G_+\{\calf_t,\calf_{t-1}\}(\Gamma)$ is a free $\Omega_*$-module.
 \item The boundary homomorphism is surjective $$\Omega^G_+\{\calf_s,\calf_t\}(\Gamma) \stackrel{\partial
}{\to} \Omega^G_-\{\calf_t,\calf_{t-1}\}(\Gamma).$$
 \end{itemize}
 
Let us see that these properties imply that $\Omega^G_*(\Gamma)$ is a free $\Omega_*$-module on even-dimensional generators.
Since $\Omega^G_*\{\calf_s,\calf_0\}(\Gamma)$ is a free $\Omega_*$-module on even-dimensional generators, the long exact sequence associated to the families of groups $\calf_{0}  \subset \calf_s$ induce the exact sequence
 $$0 \to \Omega^G_+\{\calf_0\}(\Gamma) \to \Omega^G_+(\Gamma) \to \Omega^G_+\{\calf_s,\calf_{0}\}(\Gamma) \stackrel{\partial}{\to} \Omega^G_-\{\calf_0\}(\Gamma) \to \Omega_-^G( \Gamma) \to 0.$$
 
 The unitary bordism group of free actions $\Omega^G_*\{\calf_0\}(\Gamma)$ is isomorphic
 to $\Omega_*(BG \times \prod_{i=1}^l BU(k_i))$ since both $EG \times B_GU(k_i)$ and $EG \times BU(k_i)$ classify
 $G$-equivariant complex vector bundles of rank $k_i$ over free $G$-spaces. The unitary bordism groups of $BU(k_i)$
 are free $\Omega_*$-modules on even-dimensional generators, and therefore by the K\"unneth theorem we have that
 $$\Omega^G_*\{\calf_0\}(\Gamma) \cong \Omega_*(BG) \underset{\Omega_*}{\otimes} \Omega_*\left(\prod_{i=1}^l BU(k_i)\right).$$
 
 Hence we have that $\Omega^G_+\{\calf_0\}(\Gamma)$ is a free $\Omega_*$-module in even degrees and
 that $\Omega^G_-\{\calf_0\}(\Gamma)$ is all $p$-torsion.
 Consider a unitary bordism class defined by the map $h:M \to \prod_{i=1}^l BU(k_i)$ and denote by $E:= E_1 \oplus \cdots \oplus E_l$ with $E_j$ the complex vector bundle that the map $\pi_j \circ h : M \to BU(k_j)$ defines. Take the ball $B^{2n+2}_{p^s}$
 of vectors in $\IC^{n+1}$ with norm less than 1 endowed with the action of $G$ given
 by multiplication by $e^{\frac{2 \pi i}{p^s}}$ and consider the $G$-equivariant
 $\prod_{i=1}^l U(k_i)$ complex bundle that the product 
 $B^{2n+2}_{p^s} \times E \to  B^{2n+2}_{p^s}\times M$
 defines.
 
 This $G$-equivariant $\prod_{i=1}^l U(k_i)$ complex bundle is classified by a $G$-equivariant map
 $$f : B^{2n+2}_{p^s}\times M \to \Gamma$$ 
 and its $G$-equivariant unitary bordism class $[ B^{2n+2}_{p^s}\times M, f]$ belongs
 to $\Omega_+^G\{\calf_s, \calf_0\}(\Gamma)$. Its boundary is $[ S^{2n+1}_{p^s}\times M, f|_{S^{2n+1}_{p^s}\times M}]$ and it belongs to $\Omega_-^G\{ \calf_0\}(\Gamma)$.  By the K\"unneth
 isomorphism described above we know that the unitary bordism classes $[ S^{2n+1}_{p^s}\times M, f|_{S^{2n+1}_{p^s}\times M}]$ generate $\Omega_-^G\{ \calf_0\}(\Gamma)$ and therefore the boundary
 homomorphism $\Omega^G_+\{\calf_s,\calf_{0}\}(\Gamma) \stackrel{\partial}{\to} \Omega^G_-\{\calf_0\}(\Gamma)$ is surjective. This implies that $\Omega_-^G(\Gamma)$ is trivial.
 
 Since $\Omega^G_+\{\calf_0\}(\Gamma) \cong \Omega_+(\prod_{i=1}^l BU(k_i))$  is a free $\Omega_*$-module,
 and by hypothesis  $\Omega^G_+\{\calf_s,\calf_{0}\}$ also, then
 it implies that $\Omega^G_+(\Gamma)$ is a free $\Omega_*$-module.
 
 In particular we have that the $G$-equivariant unitary bordism group $\Omega_*^G$ is a free $\Omega_*$-module on even-dimensional generators.
 
 Now let us sketch the proof of the properties cited above. Let us assume that the properties hold for cyclic groups
 of order less than $p^s$ and let us proceed by induction on the families of subgroups of $G$. For the adjacent pair of families
 $(\calf_s, \calf_{s-1})$ differing by the group $G$, we know by Theorem \ref{theorem adjacent} that
    $\Omega_*^G \{\calf_s,\calf_{s-1}\}(\Gamma)$ is a direct sum of groups $\Omega_*(\Gamma^G \times \Gamma')$
    where both $\Gamma^G$ and $\Gamma'$ are products of classifying spaces of unitary groups. Therefore
    $\Omega_*^G  \{\calf_s,\calf_{s-1}\}(\Gamma)$ is a free $\Omega_*$-module on even-dimensional generators and we have 
    started our induction.

    Now let us assume that the properties hold for the pair of families $(\calf_s, \calf_j)$ for $s >j \geq t$. Therefore we get the following exact sequence of groups
   \begin{align} \nonumber
   0 \to \Omega^G_+\{\calf_t, \calf_{t-1}\}(\Gamma)  \to \Omega^G_+&\{\calf_s, \calf_{t-1}\}(\Gamma) \to  \Omega^G_+\{\calf_s,\calf_t\}(\Gamma) \stackrel{\partial}{\to} \\ &\Omega^G_-\{\calf_t, \calf_{t-1}\}(\Gamma) \to \Omega_-^G\{\calf_s, \calf_{t-1}\}( \Gamma) \to 0. \label{LES t t-1}
   \end{align}

     Since the pair of families $(\calf_t, \calf_{t-1})$ differ by the cyclic group $H$ of order $p^t$, $G/H$ is
     a cyclic group of order $p^{s-t}$, and      $\Gamma^{H}$ is a product of classifying spaces of the form
     $B_{G/H}U(k)$, then by Theorem \ref{theorem adjacent}  there is an isomorphism 
  \begin{align} \label{iso t t-1}
     \Omega^G_*\{\calf_t, \calf_{t-1}\}(\Gamma) \cong \bigoplus_{k \geq 0}\Omega_{*-2k}^{G/H}\{\calf_0\} (\Gamma^{H} \times C_{G,H}(k))
     \end{align}
     where both $\Gamma^{H}$ and $C_{G,H}(k)$ are disjoint unions of products of classifying spaces of the form $B_{G/H}U(k)$.

     Therefore we know that  $\Omega^G_+\{\calf_t, \calf_{t-1}\}(\Gamma)$ is a free $\Omega_*$-module
     and by the induction hypothesis  we know that the boundary map 
     \begin{align} \label{induction sujectivity}
     \Omega_+^{G/H}\{\calf_{s-t},\calf_0\} (\Gamma^{H} \times C_{G,H}(k)) \stackrel{\partial}{\to} \Omega_-^{G/H}\{\calf_0\} (\Gamma^{H} \times C_{G,H}(k))
     \end{align}
     is surjective.  A bordism class in $\Omega^G_-\{\calf_t, \calf_{t-1}\}(\Gamma)$ can be represented
     by a class $[D(E), f: D(E) \to \Gamma]$ where $D(E)$ is the disk bundle  of a $G$-equivariant vector bundle $E \to M$ over a manifold $M$ on which $H$ acts trivially
     and $G/H$ acts freely, and such that the trivial representation of $H$ does not appear on the fibers of $E$. 
    This bundle is classified by a $G/H$-equivariant map $h: M \to C_{G,H}(k)$ for some $k$, and the bordism
    class $[M, f|_M \times h : M \to \Gamma^{H}\times C_{G,H}(k)]$  lives in 
    $\Omega_-^{G/H}\{\calf_0\} (\Gamma^{H} \times C_{G,H}(k))$. By the surjectivity
   of \eqref{induction sujectivity} there is a bordism class $[Z, F \times \tilde{h}: Z \to  \Gamma^{H}\times C_{G,H}(k)]$
   in $ \Omega_+^{G/H}\{\calf_{s-t},\calf_0\} (\Gamma^{H} \times C_{G,H}(k))$ such that
   $\partial Z = M$, $F|_{M}=f|_M$ and $\tilde{h}|_M=h$. Let $p: V \to Z$ denote the $G$-equivariant vector bundle
   over $Z$ that the map $\tilde{h}$ defines and note that the bordism class $[D(V), F \circ p : D(V) \to \Gamma]$ 
    defines an element in $\Omega_*^G\{\calf_s, \calf_{t}\}(\Gamma)$ since the trivial $H$-representation
    does not appear on the fibers of $V$ and the action of $G/H$ over $M$ is free.
    The boundary of $D(V)$ is the union of the sphere bundle $S(V)$ and $D(V)|_{M}= D(E)$, but since
    $S(V)$ is $\calf_{t-1}$-free we have that 
    $$ \partial [D(V), F \circ p : D(V) \to \Gamma] = [D(E), p|_{E} \circ f|_{M}: D(E) \to \Gamma]=[D(E), f: D(E) \to \Gamma]$$
and therefore the boundary map
$$\Omega_+^{G}\{\calf_s, \calf_t\}(\Gamma) \stackrel{\partial}{\to} \Omega_-^{G}\{\calf_t, \calf_{t-1}\}(\Gamma)$$
is surjective.

Now, the group $\Omega^G_*\{\calf_t, \calf_{t-1}\}(\Gamma)$ has projective dimension 1 as a $\Omega_*$-module.
This follows from the following two facts, first that $\Omega_*(B(G/H))$ has projective dimension 1 over $\Omega_*$, and second that equation \eqref{iso t t-1} induces the isomorphism
$$\Omega^G_*\{\calf_t, \calf_{t-1}\}(\Gamma) \cong \bigoplus_{k \geq 0}\Omega_{*-2k}(B(G/H) \times Z)$$
where $Z$ is a product of copies of $BU(k)$'s. By Schanuel's lemma we know that the kernel of the boundary map $\partial$
is projective and therefore free \cite[Proposition 3.2]{ConnerSmith}.
By the long exact sequence described in \eqref{LES t t-1} we deduce that $\Omega_*^{G}\{\calf_s, \calf_{t-1} \}(\Gamma)$
is a free $\Omega_{*}$-module on even-dimensional generators  and we conclude that the properties also hold for the pair of families
$(\calf_s, \calf_{t-1})$. 

 Therefore the evenness conjecture holds for cyclic $p$-groups \cite[Theorem 1']{Landweber-cyclic}.

\subsection{General abelian $p$-groups} The argument to show the evenness conjecture for general abelian $p$-groups
is more elaborate than the one done above for cyclic $p$-groups. We will follow the original proof of Stong done in \cite{Stong-complex} in which the author uses very cleverly the Thom isomorphism and the long
exact sequence for pairs of spaces in order to understand the long exact sequence for a pair of families once
restricted to a special kind of actions on manifolds. Here we shorten the original proof and we
highlight its main ingredients. 

  Let $G = H \times \IZ/q$ with $q=p^s$ such that all elements in $H$ have order less or equal than $p^s$ and
  let $$\Gamma:=  \prod_{i=1}^lB_GU(k_i)$$ be a product of spaces of the form $B_GU(k)$. We will
  show by induction on the order of the group $G$ that the bordism group $\Omega_*^G(\Gamma)$ is a 
  free $\Omega_*$-module on even-dimensional generators. Therefore let us assume that
 $\Omega_*^K(\Gamma')$ is a free $\Omega_*$-module on even-dimensional generators
 for all $p$-groups of order less than the order of $G$ and $\Gamma'$ any product of
 classifying spaces of the form $B_{K}U(l)$.

  Following the notation of Stong in \cite{Stong-complex} let us consider the following families of
  subgroups of $G$:
  \begin{itemize}
  \item $\calf_a$ is the family of all subgroups of $G$, 
  \item $\calf_s$ is the family of subgroups whose intersection with $\IZ/q$ is a proper subgroup of $W$, i.e. $\calf_s: = \{ W \subset H \times \IZ/q \, | \, \{1\} \times \IZ/q \not\subset W \}$,
  \item $\calf_f$ is the family of subgroups whose intersection with $\IZ/q$ is the unit subgroup, i.e. $\calf_f: = \left\{ W \subset H \times \IZ/q \, | \, \left( \{1\} \times \IZ/q \right) \cap W = \{(1,1)\} \right\}$.
  \end{itemize}

  A manifold $M$ is $\calf_s$-free if  for every $x \in M$ the isotropy group $(\IZ/q)_x \neq \IZ/q$, and
  it is $\calf_f$-free if the restriction of the action to $\IZ/q$ is free.
  
   The classifying space $E \calf_f$ has a free $\IZ/q$-action and can be understood
  as the universal $H$-equivariant $\IZ/q$-principal bundle $E_H \IZ/q$  \cite[Theorem 11.4]{LueckUribe} . 
  Hence $E \calf_f = E_H \IZ/q$ and its quotient $E \calf_f / (\IZ/q) = B_H \IZ/q$ is the
  classifying space of $H$-equivariant $\IZ/q$-principal bundles. By the isomorphism of
  \eqref{iso geometric families EF}, and since the action of $\IZ/q$ is free, we get the following isomorphisms
  \cite[Proposition 3.1]{Stong-complex}:  
   \begin{align}\label{iso gf}
   \Omega_*^G\{ \calf_f\}(X) \cong \Omega_*^G(X \times E_H \IZ/q) \cong \Omega_*^H(X \times_{\IZ/q} E_H \IZ/q).
   \end{align}
  
  Since both spaces $E_H \IZ/q \times B_GU(k_i)$ and $E_H \IZ/q \times B_H U(k_i)$
  classify $H \times \IZ/q$-equivariant $U(k_i)$-principal bundles over spaces with free $\IZ/q$-action,
  we may take the maps
  $$B_H U(k_i) \to B_G U(k_i) \to B_H U(k_i),$$
  where the left hand side map classifies the $G$-equivariant complex bundles such that
  the action of $\IZ/q$ is trivial over the total space of the bundle, and the right hand side is the one that
  forgets the $\IZ/q$-action, thus producing $G$-equivariant homotopy equivalences
  $$E_H \IZ/q \times E_H U(k_i) \stackrel{\simeq}{\to} E_H \IZ/q \times E_G U(k_i) \stackrel{\simeq}{\to} E_H \IZ/q \times E_H U(k_i).$$
  
  If we denote by $$\Gamma': = \prod_{i=1}^lB_HU(k_i)$$ and the map $\iota: \Gamma' \to \Gamma$ is the one
  that classifies trivial $\IZ/q$-bundles over $H$-spaces, then the argument above implies that the isomorphism  \cite[Proposition 3.2]{Stong-complex} holds:
  \begin{align}\label{iso gf for classifying space}
   \Omega_*^H\left( B_H \IZ/q \times \Gamma' \right) \cong  \Omega_*^G\left( E_H \IZ/q \times \Gamma' \right) \underset{\cong}{\stackrel{\iota_*}{\to}}\Omega_*^G\{ \calf_f\}\left(\Gamma \right).
   \end{align}
  
 Let $T$ be the generator of the group $\IZ/q$ and denote by $\IZ/p^t$ the subgroup generated by $T^{p^{s-t}}$.
 A manifold is $\calf_f$-free if and only if $T^{p^{s-1}}$ acts freely  and therefore a $(\calf_a, \calf_f)$-manifold $M$
 can be localized to the normal bundle of the fixed point set $M^{\IZ/p}$ of the subgroup $\IZ/p$.
 The normal bundle is a $G$-equivariant complex bundle over the trivial $\IZ/p$ space and once it is classified
 to the appropriate spaces $C_{G, \IZ/p}(k)$ of \eqref{space CGA(k)} we obtain the isomorphism  \cite[Proposition 3.4]{Stong-complex}:    
  \begin{align}\label{iso g,gf}
   \Omega_*^G\{\calf_a, \calf_f\}(X) \cong \bigoplus_{0\leq 2k \leq *} \Omega_{*-2k}^{G/ (\IZ/p)}\left( X^{\IZ/p} \times C_{G,  \IZ/p}(k) \right).
   \end{align}
  
  Applying the previous isomorphism to $\Gamma = \prod_{i=1}^lB_GU(k_i)$ we obtain that 
  $$ \Omega_*^G\{\calf_a, \calf_f\}(\Gamma)$$ is a free $\Omega_*$-module on even-dimensional generators
  since both $\Gamma^{\IZ/p}$ and $C_{G,  \IZ/p}(k)$ are products of spaces of the form $B_{G/ ( \IZ/p)} U(l)$
  and by induction we assumed that the evenness conjecture was true for groups of order less than the one of $G$
  and spaces of this type. Therefore the long exact sequence for the pair of families $(\calf_a, \calf_f)$ becomes:
  \begin{align} \label{exact sequence g,gf}
  0 \to \Omega_+^G\{\calf_f\}(\Gamma) \to \Omega_+^G(\Gamma){\to} \Omega_+^G\{\calf_a,\calf_f\}(\Gamma)
   \stackrel{\partial}{\to} \Omega_-^G\{\calf_f\}(\Gamma) \to \Omega_-^G(\Gamma) \to 0.
  \end{align}

A $(\calf_a, \calf_s)$-free manifold $M$ once restricted to the action of $\IZ/q$ becomes
a $\IZ/q$-manifold on which the boundary has no fixed points of the whole group. Therefore
the manifold can be localized on the normal bundle of the fixed point set $M^{\IZ/q}$ and the information
of the normal bundle can be recorded with appropriate maps to the classifying spaces
$C_{G, \IZ/q}(k)$ of \eqref{space CGA(k)}. Following the same proof as in Theorem \ref{theorem adjacent}
one obtains the following isomorphism \cite[Proposition 3.3]{Stong-complex}:
 \begin{align}\label{iso g,gs}
   \Omega_*^G\{\calf_a, \calf_s\}(X) \cong \bigoplus_{0\leq 2k \leq *} \Omega_{*-2k}^{H}\left( X^{\IZ/q} \times C_{G,  \IZ/q}(k) \right).
   \end{align}
   Since both $\Gamma^{\IZ/q}$ and $C_{G,  \IZ/q}(k)$ are products of spaces of the form $B_HU(l)$, by 
   the induction hypothesis we obtain that $ \Omega_*^G\{\calf_a, \calf_s\}(\Gamma)$ is a free $\Omega_*$-module
   on even-dimensional generators. 

In order to understand the image of the boundary map of \eqref{exact sequence g,gf} Stong restricted 
the equivariant bordism groups to manifolds with a special type of $G$ action. Stong noticed that
the image of the boundary map could be determined by restricting to manifolds on which the $\IZ/q$-fixed
points are of codimension 2 and therefore he studied the class of {\it{special $G$ actions}}.

\begin{definition}
Let $G= H \times \IZ/q$ be a finite abelian group. The class of {\it{special $G$ actions}} is the collection of
$G$-equivariant unitary manifolds $M$ satisfying:
\begin{itemize}
\item The restriction to a $\IZ/q$-action is semi-free, i.e. for each $x \in M$ the isotropy group $(\IZ/q)_x$
is either $\IZ/q$ or $\{1\}$.
\item The set $M^{\IZ/q}$ of fixed point sets has codimension 2 in $M$  and $\IZ/q$ acts on the normal bundle
of $M^{\IZ/q}$ so that the generator $T$ of $\IZ/q$ acts by multiplication by $e^{\frac{2 \pi i }{q}}$, or the fixed
point set $M^{\IZ/q}$ is empty.
\end{itemize}
\end{definition}

The class of special $G$ actions is sufficiently large to permit all constructions done in section \ref{section bordism}, and for
a pair of families $(\calf, \calf')$ in $G$ we denote by $\overline{\Omega}_*^G\{\calf,\calf'\}$ the equivariant
homology theory defined by using only special $G$ actions. The inclusion of special $G$ actions in the class of all $G$
actions defines natural transformations of homology theories
$$I_* \colon \overline{\Omega}_*^G\{\calf,\calf'\} \to {\Omega}_*^G\{\calf,\calf'\}$$
preserving the relations between these functors. The $G$-equivariant unitary bordism groups of special $G$ actions satisfy
the following properties:

\begin{itemize}
\item[(i)] The natural transformation $$I_* \colon \overline{\Omega}_*^G\{\calf_f\} \stackrel{\cong}{\to} {\Omega}_*^G\{\calf_f\}$$is an equivalence since every $\calf_f$ action  is a special $G$ action.
\item[(ii)] The inclusion $(\calf_a, \calf_f) \subset (\calf_a, \calf_s)$ induces an isomorphism 
$$\overline{\Omega}_*^G\{\calf_a, \calf_f\}(X) \stackrel{\cong}{\to} \overline{\Omega}_*^G\{\calf_a, \calf_s\}(X)$$
since $\calf_s$-free special $G$ actions are $\calf_f$-free.
\item[(iii)] From the equation \eqref{iso g,gs} we get the isomorphism
$$\overline{\Omega}_*^G\{\calf_a, \calf_s\}(X) \cong \Omega_{*-2}^H(X^{\IZ/q} \times B_HU(1)),$$
thus implying that $\overline{\Omega}_*^G\{\calf_a, \calf_s\}(X)$ maps isomorphically to a direct summand in ${\Omega}_*^G\{\calf_a, \calf_s\}(X)$.
\item[(iv)] For $\Gamma:=\prod_{i=1}^lB_GU(k_i)$ the induction hypothesis implies that
$\overline{\Omega}_*^G\{\calf_a, \calf_f\}(\Gamma)$ is a free $\Omega_*$-module on even-dimensional generators.
Therefore the canonical maps
$$\overline{\Omega}_*^G\{\calf_a, \calf_f\}(\Gamma) \to {\Omega}_*^G\{\calf_a, \calf_f\}(\Gamma) \to {\Omega}_*^G\{\calf_a, \calf_s\}(\Gamma) $$ imply that $\overline{\Omega}_*^G\{\calf_a, \calf_f\}(\Gamma)$ also maps isomorphically to a direct summand in ${\Omega}_*^G\{\calf_a, \calf_f\}(\Gamma)$.
\end{itemize}

Let us now concentrate on understanding the five-term exact sequence restricted to special $G$ actions
 \begin{align} \label{exact sequence special g,gf}
  0 \to \overline{\Omega}_+^G\{\calf_f\}(\Gamma) \to \overline{\Omega}_+^G(\Gamma){\to} \overline{\Omega}_+^G\{\calf_a,\calf_f\}(\Gamma)
   \stackrel{\partial}{\to} \overline{\Omega}_-^G\{\calf_f\}(\Gamma) \to \overline{\Omega}_-^G(\Gamma) \to 0.
  \end{align}
Note that the map $\iota_*:\Gamma' \to \Gamma$ induces the commutative diagram
$$
\xymatrixrowsep{0.5cm}
\xymatrix{
\Omega_{*-2}^H(\Gamma^{\IZ/q} \times B_HU(1)) \ar[r]^{\cong} & \overline{\Omega}_*^G\{\calf_a,\calf_f\}(\Gamma) \ar[r]^\partial & \ \overline{\Omega}_{*-1}^G\{\calf_f\}(\Gamma) \\
\Omega_{*-2}^H(\Gamma' \times B_HU(1)) \ar[r]^{\cong} \ar[u]&  \overline{\Omega}_*^G\{\calf_a,\calf_f\}(\Gamma') \ar[u]^{\iota_*} \ar[r]^\partial & \ \overline{\Omega}_{*-1}^G\{\calf_f\}(\Gamma') \ar[u]^{\iota_*}_{\cong}
}$$
where the middle homomorphism $\iota_*$ maps isomorphically $\overline{\Omega}_*^G\{\calf_a,\calf_f\}(\Gamma')$ 
into a direct summand in $\overline{\Omega}_*^G\{\calf_a,\calf_f\}(\Gamma)$ since $\Gamma'$ is mapped
to one connected component of the fixed point set $\Gamma^{\IZ/q}$. Therefore the image of the boundary
homomorphism $\partial$ is the same in both cases. 

In what follows we will study the induced boundary homomorphism 
\begin{align} \label{boundary homomorpshim gamma'}
\Omega_{*-2}^H(\Gamma' \times B_HU(1)) \to  \overline{\Omega}_{*-1}^G\{\calf_f\}(\Gamma') \cong 
{\Omega}_{*-1}^H(\Gamma' \times B_H \IZ/q)
\end{align} using the Thom isomorphism, the long exact sequence for pairs and 
a particular model for $E_H\IZ/q$.

Let ${\bf{C}}^\infty_H$ be a countable direct sum of all complex irreducible $H$-representations and consider
the $\IZ/q$ action on ${\bf{C}}^\infty_H$ such that the generator $T$ of $\IZ/q$ acts by mutliplication
with $e^{\frac{2 \pi i}{q}}$. The sphere $S({\bf{C}}^\infty_H)$ of vectors of norm 1 is an $G=\IZ/q \times H$
space on which $\IZ/q$ acts freely and moreover is a $\calf_f$-free space. Since the non-empty fixed point
sets are infinite dimensional spheres we know that this sphere $S({\bf{C}}^\infty_H)$ is a model for
$E_H\IZ/q$. The Grassmannian $Gr_1 {\bf{C}}^\infty_H$ is a model for $B_HU(1)$ since
$\IZ/q$ acts trivially on the one dimensional vector spaces, and $\IZ/q$ acts on the fibers of the canonical line bundle
$\gamma_HU(1) \to B_HU(1)$ by multiplication by $e^{\frac{2 \pi i}{q}}$. 
To simplify the notation denote by $$\gamma_1 : = \gamma_HU(1)$$  and note that $S({\bf{C}}^\infty_H) \cong S(\gamma_1)$ where $S(\gamma_1)$ denotes the sphere bundle of $\gamma_1$.

Consider now the line bundle $\gamma_1^{\otimes q}$ over $B_HU(1)$ which is the $q$-fold tensor product of $\gamma_1$.
The diagonal map $$ \Delta: \gamma_1 \to \gamma_1^{\otimes q}, \ \ v \mapsto v \otimes \cdots \otimes v$$
is a $q$ to 1 map on the fibers of the line bundles and therefore it induces an $H$-equivariant homeomorphism
$$S(\gamma_1)/( \IZ/q) \cong S(\gamma_1^{\otimes q}).$$
Therefore we may take either $S(\gamma_1)/( \IZ/q)$ or $S(\gamma_1^{\otimes q})$ as a model
for $B_H\IZ/q$. The Thom isomorphism 
$$\Omega_*^H(( D(\gamma_1^{\otimes q}),  S(\gamma_1^{\otimes q})) \times \Gamma') \cong 
\Omega_{*-2}^H(B_HU(1) \times \Gamma')$$
 together with the long exact sequence for the pair $( D(\gamma_1^{\otimes q}),  S(\gamma_1^{\otimes q}))$
 and the induction hypothesis provides  a four-term exact sequence 
 \begin{align*}
 0 \to \Omega_+^H( \Gamma' \times S(\gamma_1^{\otimes q})) \to  \Omega_+^H&( \Gamma' \times B_HU(1)) \to  \\
& \Omega_{+}^H( \Gamma' \times B_HU(1)) \to \Omega_{-}^H( \Gamma' \times S(\gamma_1^{\otimes q})) \to 0,\end{align*}
where the right hand side homomorphism is precisely the one of
 \eqref{boundary homomorpshim gamma'}. Therefore we obtain that the boundary homomorphism
 of \eqref{boundary homomorpshim gamma'} is surjective, and since by the induction hypothesis
 $ \Omega_+^H( \Gamma' \times B_HU(1))$ is a free $\Omega_*$-module, we conclude that $ \Omega_+^H( \Gamma' \times S(\gamma_1^{\otimes q}))$ is also a free $\Omega_*$-module.
 
 Therefore we obtain the following commutative diagram with exact rows
 $$
 \xymatrixrowsep{0.4cm}
\xymatrixcolsep{0.4cm}
 \xymatrix{
0 \ar[r] & {\Omega}_{+}^G\{\calf_f\}(\Gamma)  \ar[r] & {\Omega}_+^G(\Gamma) \ar[r] & {\Omega}_+^G\{\calf_a,\calf_f\}(\Gamma) \ar[r]^\partial &  {\Omega}_{-}^G\{\calf_f\}(\Gamma) \ar[r] &0 \\
0 \ar[r] &\overline{\Omega}_{+}^G\{\calf_f\}(\Gamma) \ar[r] \ar[u]_{\cong}&\overline{\Omega}_+^G(\Gamma) \ar[r] \ar[u]& \overline{\Omega}_+^G\{\calf_a,\calf_f\}(\Gamma) \ar[r]^\partial \ar@{^{(}->}[u]&  \overline{\Omega}_{-}^G\{\calf_f\}(\Gamma) \ar[u] _\cong \ar[r] &0\\
0 \ar[r] & \overline{\Omega}_{+}^G\{\calf_f\}(\Gamma') \ar[r] \ar[u]^{\iota_*}_{\cong} &\overline{\Omega}_+^G(\Gamma') \ar[r] \ar[u]&  \overline{\Omega}_+^G\{\calf_a,\calf_f\}(\Gamma') \ar@{^{(}->}[u] \ar[r]^\partial & \ \overline{\Omega}_{-}^G\{\calf_f\}(\Gamma') \ar[u]^{\iota_*}_{\cong} \ar[r] &0,
}$$
 thus implying that $\Omega_-^G(\Gamma)=0$ and that $\Omega_+^G(\Gamma)$ is a free $\Omega_*$-module
 since both ${\Omega}_{+}^G\{\calf_f\}(\Gamma)$ and ${\Omega}_+^G\{\calf_a,\calf_f\}(\Gamma)$ are free
 $\Omega_*$-modules.
 
 Therefore the evenness conjecture holds for finite abelian $p$-groups \cite{Stong-complex}.
       
       \subsection{The general case} The proof of the evenness conjecture for general finite abelian groups
       was done by Ossa \cite{Ossa} and is based on the proof of Stong for $p$-groups and 
       appropriate localizations at different primes. For a finite abelian group $K$ denote by $Z_K:= \IZ[1/|K|]$ the localization of the integers at the ideal generated by the order of $K$.

Let $G = K \times L$ with $K$ and $L$ finite abelian with $|K|$ and $|L|$ relatively prime and consider
the homomorphism $\Omega_*^{K \times L}\{\calf\} \to \Omega_*^{ L}\{\calf\}$ which forgets the $K$ action
and $\calf$ is any family of subgroups of $L$ . Let us 
show that the localized homomorphism $$\Omega_*^{K \times L}\{\calf\}(\Gamma)\otimes Z_K \to \Omega_*^{ L}\{\calf\}(\Gamma) \otimes Z_K$$ is an isomorphism whenever $\Gamma:= \prod_{i=1}^l B_GU(k_i)$. Let us proceed
by induction over $L$ and over the family $\{\calf\}$. 

For the trivial family $\calf=\{\{1\}\}$ we obtain the isomorphism
$$\Omega_*(BK \times BL \times \prod_i BU(k_i)) \otimes Z_K \stackrel{\cong}{\to}
\Omega_*( BL \times \prod_i BU(k_i))\otimes Z_K$$
since $\Omega_*(BK) \otimes Z_K \cong \Omega_* \otimes Z_K$.

Whenever the adjacent pair of families $( \calf, \calf')$ differ by $H \subset L$ we obtain the homomorphism
of long exact sequences
   $$
 \xymatrixrowsep{0.4cm}
\xymatrixcolsep{0.4cm}
 \xymatrix{
 {} \ar[r] & \Omega_*^{K \times L}\{\calf'\}(\Gamma) \ar[r] \ar[d] & \Omega_*^{K \times L}\{\calf\}(\Gamma) \ar[r] \ar[d] &
  \Omega_*^{K \times L/H}\{\{1\}\}(\Gamma^H \times \Gamma') \ar[r]\ar[d]  & \\
  {} \ar[r] & \Omega_*^{ L}\{\calf'\}(\Gamma) \ar[r] & \Omega_*^{ L}\{\calf\}(\Gamma) \ar[r] &
  \Omega_*^{ L/H}\{\{1\}\}(\Gamma^H\times \Gamma') \ar[r] & 
 }  $$
 with $\Gamma'$ a disjoint union of products of spaces of the form $B_{K \times L/H}U(l)$.
Tensoring with $Z_K$ induces an isomorphism on the left vertical arrow by the induction
 hypothesis on the families and an isomorphism on the right vertical arrow by the induction hypothesis on the group $L/H$.
The 5-lemma implies the desired isomorphism.

Now let  $\calf$ be any family of subgroups of $K$ and denote by $\calf \times \Phi$ the family
of subgroups of $G$ whose elements are groups $J\times H$ with $J \in \calf$ and $H$ any subgroup of $L$.
Let us show by induction on $\calf$ and on the group $K$ that the localized module
$$\Omega_*^{K \times L}\{\calf \times \Phi\}(\Gamma)\otimes Z_K $$
is a free $\Omega_* \otimes Z_K$-module on even-dimensional generators. Whenever $\calf$ is the trivial family we have shown above
that $$\Omega_*^{K \times L}\{\{1\} \times \Phi\}(\Gamma)\otimes Z_K \stackrel{\cong}{\to} \Omega^L(\Gamma) \otimes Z_K$$
is an isomorphism, and by the induction hypothesis we know that $\Omega^L(\Gamma)$ is a free $\Omega_*$-module on
even-dimensional generators.

If the adjacent pair of families $(\calf, \calf')$ differ by the subgroup $J$, then we obtain the long exact sequence  
$$ \cdots \to \Omega_*^{K \times L}\{\calf'\times \Phi\}(\Gamma)  \to  \Omega_*^{K \times L}\{\calf \times \Phi\}(\Gamma)  \to 
 \Omega_*^{K/J}\{\{1\}\}(\Gamma^{J\times L} \times \Gamma'')  \to \cdots$$
where $\Gamma''$ is a disjoint union of spaces of the form $B_{K/J}U(l)$. Tensoring with $Z_K$ gives us free $\Omega_* \otimes Z_K$-modules on even-dimensional generators on the left hand side by the induction on families and free $\Omega_* \otimes Z_K$-modules on even-dimensional generators on the right hand side by the induction on the group $K$. By \cite[Proposition 3.2]{ConnerSmith} projective $\Omega_* \otimes Z_K$-modules
are free, hence
the  middle term is also a free $\Omega_* \otimes Z_K$-module on even-dimensional generators.

Therefore we have proved that if $\Omega^L(\Gamma)$ is a free $\Omega_*$-module then $\Omega^{K \times L}(\Gamma)
\otimes Z_K$ is a free $\Omega_* \otimes Z_K$-module. Let us now write $G = P_1 \times \cdots \times P_k$
with $P_i$ its Sylow $p_i$-subgroup. Since the evenness conjecture holds for $p$-groups, we have that
$\Omega_*^{P_i}(\Gamma)$ is a free $\Omega_*$-module and therefore $\Omega_*^G(\Gamma) \otimes \IZ[1/[G :P_i]]$ is a free $\Omega_* \otimes \IZ[1/[G :P_i]]$-module. Since the numbers $[G :P_i]$ are relatively prime
it follows that $\Omega_*^G(\Gamma)$ is a free $\Omega_*$-module.
              
              Therefore the evenness conjecture holds for finite abelian groups \cite{Ossa}.
              
       \section{The equivariant unitary bordism groups for non-abelian groups.} \label{section non abelian}
           
The evenness conjecture has been shown to be true for the dihedral groups $D_{2p}$ with $p$ prime
by \'Angel, G\'omez and the author \cite{AngelGomezUribe}, for groups of order $pq$ where $p$ and $q$ are different primes by Lazarov
\cite{Lazarov} and for the more general case of groups for which all its Sylow subgroups are cyclic by Rowlett \cite{Rowlett-metacyclic}.
In these cases the group $G$ is a semidirect product $\IZ/r \rtimes \IZ/s$ of cyclic groups with $r$ and $s$ relatively prime
\cite[Theorem 9.4.3]{Hall},
 and the study of the equivariant unitary bordism groups is also carried out by calculating the equivariant
unitary bordism groups of adjacent pairs of families of subgroups as is done in the cyclic group case of section \ref{subsection cyclic}. 

The main tool used by Rowlett to study the case on which all Sylow subgroups are cyclic is the equivariant unitary spectral sequence
constructed by himself in \cite[Prop.osition 2.1]{Rowlett}. Suppose that $A$ is a normal subgroup of $G$ and that $Q=G/A$.
A family $\calf$ of subgroups of $A$ is called $G$-invariant if it is closed under conjugation
by elements of $G$.  Consider a pair $(\calf, \calf')$ of $G$-invariant families of $A$ and note that 
$\Omega_*^A\{\calf, \calf'\}$ becomes a $Q$-module in the following way. Consider an $A$-manifold $M$
with action $\theta: A \times M \to M$ and  take an element $g \in G$. Define a new action on $M$
by the map $g_* \theta: A \times M \to M$, $ g_*(a,m) := \theta(g^{-1}ag,m)$ and denote the action
of $g$ on the  bordism class $[M, \theta]$ by $\overline{g}[M,\theta]: =[M, g_* \theta]$. This action
is trivial on elements of $A$ and therefore it boils down to an action of $Q$.
Then there is a first quadrant spectral sequence $E^r$ converging to $\Omega_*^G\{\calf, \calf'\}$
whose second page is 
$$E^2_{p,q}\cong H_p(Q,\Omega_q^A\{\calf, \calf'\}).$$

In the case that both groups $A$ and $Q$ are cyclic of relatively prime order, the action of $Q$ on 
$\Omega_+^A\{\calf, \calf'\}$ factors through a permutation action on the free generators and therefore the 
second page is not difficult to calculate.
If we take the family $\calf_A$ of all subgroups of $A$, the second page of the spectral sequence
becomes $H_q(Q, \Omega_q^A)$, and since $\Omega_*^A$ is a free $\Omega$-module on even-dimensional
generators, then we obtain that $\Omega_+^G\{\calf_A\}$ is a free $\Omega_*$-module. Moreover, the same explicit construction
carried out in section \ref{subsection cyclic} can be adapted in this case to show that the long exact sequence
associated to the pair of families $\{\calf_a, \calf_A\}$, with $\calf_a$ the family of all subgroups, becomes
$$ 0 \to \Omega_+^G\{\calf_A \} \to \Omega_+^G \to \Omega_+^G\{\calf_a, \calf_A\} \stackrel{\partial}{\to} \Omega_-^G\{\calf_A\} \to 0.$$

The same argument as in \eqref{iso g,gf} shows that  $\Omega_*^G\{\calf_a, \calf_A\}$ is a free $\Omega_*$-module
on even-dimensional generators and therefore we conclude that $\Omega_-^G$ is zero and $\Omega_+^G$ is a free $\Omega_*$-module.

The spectral sequence defined above can also be used in order to understand the torsion-free
part of the $G$-equivariant unitary bordism groups for any abelian group. Take any subgroup $A$ of $G$
and let $(\calf_A, \calf_A')$ be the adjacent pair of families of $G$ which differ by the conjugacy class
of $A$. Tensoring with the rationals we obtain an isomorphism
$$\Omega_*^G \{\calf_A, \calf_A'\} \otimes \IQ \cong \Omega_*^A \{\calf_A, \calf_A'\}^{W_A} \otimes \IQ$$
where the right hand side consists of the $W_A$-invariant part and $W_A:= N_A/A$. Since 
$\Omega_*^A \{\calf_A, \calf_A'\}$ is a free $\Omega_*$-module on even-dimensional generators
we obtain the isomorphism
$$\Omega_*^G \otimes \IQ \cong \bigoplus_{(A)}\Omega_*^A \{\calf_A, \calf_A'\}^{W_A} \otimes \IQ$$
where $(A)$ runs over the conjugacy classes of subgroups of $G$ \cite[Theorem 1.1]{Rowlett}, cf. \cite[Theorem 1]{tomDieck-Mackey}.
In particular the torsion-free component of $\Omega_*^G$ is all of even degree.

Apart from the non-abelian groups in which all their Sylow subgroups are cyclic, there is no other finite non-abelian group
for which the evenness conjecture has been shown to hold.

 The main difficulty lies in the understanding of the equivariant bordism groups $\Omega_*^G\{\calf\}( {}^{\widetilde{G}}B_GU(n))$ of the classifying spaces ${}^{\widetilde{G}}B_GU(n)$ associated
 to $\IS^1$-central extensions $\widetilde{G}$ of $G$ for different families $\calf$ of subgroups. These
 bordism groups are the ones appearing once we try to calculate the equivariant unitary bordism groups
 for adjacent pair of families. Any development in the understanding of these 
 equivariant unitary bordism groups will shed light on the proof of the evenness conjecture for a bigger class
 of groups. 
 
 \section{Conclusion} \label{section conclusion}
 The evenness conjecture for equivariant unitary bordism has been an important question in algebraic topology
 for more than forty years. The conjecture has been proved to hold only for compact abelian Lie groups and finite
 groups for which all their Sylow subgroups are cyclic, for all other groups the conjecture remains open. 
 We do hope that the present summary of known results will help settle the conjecture in the near future.

\bibliographystyle{abbrv} 
  \bibliography{Evenness-bibliography}
\end{document}